\documentclass{article}
\usepackage{latexsym, amsxtra, amsmath, amssymb, amsthm, amscd}

\usepackage{amssymb}
\usepackage{indentfirst}
\usepackage{amscd}
\newtheorem{theorem}{Theorem}[section]
\newtheorem{lem}[theorem]{Lemma}
\newtheorem{cor}[theorem]{Corollary}

\def\Ass{\text{Ass}}

\def\gr{\text{gr}}
\def\depth{\text{depth}}
\def\H{\text{H}}
\def\m{\frak m}
\def\p{\frak p}
\def\q{\frak q}
\def\qa{\frak q(\alpha)}
\def\x{\underline x}

\def\D{\mathcal D}
\def\pds{\bigcap\limits_{\alpha\in\Lambda_{s,n}}(y_1^{\alpha_1},\ldots,y_s^{\alpha_s})M}
\def\pdm{\bigcap\limits_{\alpha\in\Lambda_{i,n}}(y_1^{\alpha_1},\ldots,y_i^{\alpha_i})M}

\def\pd{\bigcap\limits_{\alpha\in\Lambda_{d,n}}\q(\alpha)M}
\def\o{\overline}
\def\:{:_M}
\begin{document}\title{Parametric Decomposition of Powers of Parameter Ideals and Sequentially Cohen-Macaulay Modules  }        
\author{Nguyen Tu Cuong \footnote{Email: ntcuong@math.ac.vn}\ and Hoang Le Truong\footnote{Email: hltruong@math.ac.vn}\\
Institute of Mathematics\\
18 Hoang Quoc Viet Road, 10307 Hanoi, Viet Nam}        
 \date{}    
\maketitle
\begin{abstract} Let $M$ be a finitely generated module of dimension $d$ over a Noetherian local ring $(R,\m)$ and  $\q $ the parameter ideal generated by a system of parameters $\x = (x_1,\ldots , x_d)$ of $M$. For each positive integer $n$, set
$$\Lambda_{d,n}=\{ \alpha =(\alpha_1,\ldots,\alpha_d)\in\Bbb{Z}^d|\alpha_i\geqslant 1, \forall 1\leqslant i\leqslant d \text{ and } \sum\limits_{i=1}^d\alpha_i=d+n-1\}$$ 
and $\qa = (x_1^{\alpha_1},\ldots,x_d^{\alpha_d})$.  Then we prove in this note that $M$ is  a sequentially Cohen-Macaulay module if and only if there exists a certain system of parameters  $\x$ such that the equality $\q^nM=\pd$
holds true for all $n$.
 As an application of this result, we can compute the Hilbert-Samuel polynomial of a sequentially Cohen-Macaulay module with respect to certain parameter ideals.
\vspace{.2cm}
\end{abstract}

\section{Introduction} Throughout this note we denote $R$  a commutative Noetherian local  ring with the maximal ideal $\m$ and $M$ is a finitely generated $R$-module with $\dim M=d$. Let $\x = x_1,\ldots,x_d$ be a system of parameters of module $M$ and  $\q=(x_1,\ldots,x_d)$ the parameter ideal of $M$ generated  by $\x$. For each integers $n\geqslant 1$, we set $$\Lambda_{d,n}=\{(\alpha_1,\ldots,\alpha_d)\in\Bbb{Z}^d|\alpha_i\geqslant 1\text{ for all }1\leqslant i\leqslant d \text{ and } \sum\limits_{i=1}^d\alpha_i=d+n-1\}.$$
Let $\qa=(x_1^{\alpha_1},\ldots,x_d^{\alpha_d})$ for all $\alpha=(\alpha_1,\ldots,\alpha_d)\in\Lambda_{d,n}$. We say that the system of parameters $\x$   has the {\it property of  parametric decomposition}, if the equality $\q^nM=\pd$ holds true for all $n\geqslant 1$. The main purpose of this note is to study the question of when a given system of parameters of $M$ has the property of  parametric decomposition. Note that Heinzer, Ratliff and Shah \cite [Theorem 2.4]{HRS} proved that an  $R$-regular sequence  always has the property of parametric decomposition. Later, Goto and Shimoda \cite [Theorem 1.1]{GS1}  showed that the converse is also true when each element of the sequence is a non-zerodivisor in $R$.  Moreover, they gave in \cite [Theorem 1.1]{GS2} a characterization of $R$ with $\dim R\geqslant 2$, in which every system of parameters of $R$ has the property of parametric decomposition. In order to generalize  this result of Goto and Shimoda, let us recall some notions  which were defined in \cite {CC}.  A filtration $\D:\ H^0_\m(M)=D_0\subset D_1\subset\ldots\subset D_t=M$ of submodules of $M$ is  said to be a {\it dimension filtration},  if $D_{i-1}$ is the largest submodule of $D_i$ with $\dim D_{i-1}<\dim D_i$ for all $i=t,t-1,\ldots,1$.  If $D_i/D_{i-1}$ is  Cohen-Macaulay  for all $i=1,\ldots,t$, $M$ is called a {\it sequentially Cohen-Macaulay module}.  A system of parameters $\x=x_1,\ldots,x_d$ of $M$ is called a {\it good system of parameters} of $M$ if  $D_i\cap (x_{d_i+1},\ldots,x_d)M=0$ for all $i=0,\ldots,t-1$ where $d_i=\dim D_i$. Now, restrict our interest in the above question to the set of all good systems of parameters of $M$. It turns out that the property of parametric decomposition of a good system of parameters can be characterized by the sequentially Cohen-Macaulayness of the  module. The following theorem is the main result of the note.
\begin{theorem}\label{0}
The following statements are equivalent:

\noindent (i) $M$ is a sequentially Cohen-Macaulay module.

\noindent (ii)  Every good system of parameters of $M$ has the property of parametric decomposition.

\noindent (iii)  There exists a good system of parameters of $M$ having the property of parametric decomposition.
\end{theorem}
As a consequence of Theorem \ref {0} we obtain a module version for  the main result of Goto-Shimoda \cite[Theorem 1.1]{GS2}.

\begin{cor}
Let $\dim M\geqslant 2$ and $H_\frak{m}^0(M)$ the $0^{th}$ local cohomology module of $M$ with respect to the maximal ideal $\frak{m}$. Then the following statements are equivalent:

(i) $M/H^0_\m(M)$ is a Cohen-Macaulay module and $\m H^0_\m(M)=0$.

(ii) Every system of parameters of $M$ has the property of parametric decomposition.
\end{cor}

 Before  giving proofs for   Theorem 1.1 and its corollary in Section 3, we need some basic facts on good systems of parameters and sequentially Cohen-Macaulay modules, which will be summarized in Section 2. In Section 4 we shall show that  the Hilbert-Samuel polynomial of a sequentially Cohen-Macaulay module $M$ with respect to a good parameter ideal (Theorem 4.3) can be effectively computed  by using Theorem 1.1 and the dimension filtration $\D$ of $M$.
\section{Preliminaries}
Throughout this paper, $R$ is a Noetherian local commutative ring with maximal ideal $\m$ and $M$ is a finitely generated $R$-module with $\dim M=d$. Let $\x=x_1,\ldots,x_d$ be a system of parameters of module $M$ and we denote by $\q$  the ideal generated by $x_1,\ldots,x_d$. For positive integers $n, s$, we set $$\Lambda_{s,n}=\{(\alpha_1,\ldots,\alpha_s)\in\Bbb{Z}^d|\alpha_i\geqslant 1\text{ for all }1\leqslant i\leqslant s \text{ and } \sum\limits_{i=1}^s\alpha_i=s+n-1\}.$$
Let $\qa=(x_1^{\alpha_1},\ldots,x_d^{\alpha_d})$ for each $\alpha=(\alpha_1,\ldots,\alpha_d)\in\Lambda_{d,n}$.  Then $\q^n \subseteq \pd$, and if the equality $\q^n = \pd$ holds true for a system of parameters $\x$ of $M$, we say that $\x$  has the {\it property of parametric decomposition}.  Recall that a filtration $\D:\ H^0_\m(M)=D_0\subset D_1\subset\ldots\subset D_t=M$ of submodules of $M$ is  said to be a {\it dimension filtration},  if $D_{i-1}$ is the largest submodule of $D_i$ with $\dim D_{i-1}<\dim D_i$ for all $i=t,t-1,\ldots,1$, and system of parameters $\x=x_1,\ldots,x_d$ of $M$ is called a {\it good system of parameters} of $M$ if  $D_i\cap (x_{d_i+1},\ldots,x_d)M=0$ for all $i=0,\ldots,t-1$ where $d_i=\dim D_i$.

Now, let us briefly give some facts on the dimension filtration and good systems of parameters (see \cite{CC}, \cite{CN}).  Because of the Noetherian property of $M$,  the dimension filtration of $M$ exists uniquely. Therefore, in the sequel we always denote  by $$\D:\H^0_\m(M)=D_0\subset D_1\subset\ldots\subset D_t=M$$ with $\dim D_i =d_i$ the dimension filtration of $M$. In this case, we also say that the dimension filtration $\D$ of $M$ has the length $t$.  Moreover, let $\bigcap_{\p\in \Ass M}N(\p)=0$ be a reduced primary
decomposition of $0$ of $M$, then $D_i=\bigcap_{\dim(R/\p)\geqslant d_{i+1}}N(\p)$. Put $N_i=\bigcap_{\dim(R/\p)\leqslant
d_i}N(\p).$ Therefore $D_i\cap N_i=0$ and $\dim (M/N_i)=d_i$. By the
Prime Avoidance  there exists a system of parameters $\x=(x_1,\ldots,x_d)$ such that $x_{d_i+1},\ldots,x_d
\in\text {Ann }(M/N_i)$. It follows that $D_i\cap (x_{d_i+1},\ldots,x_d)M\subseteq N_i\cap D_i=0$ 
for all $i=0,\ldots , t-1$. Thus $\x=x_1,\ldots,x_d$ is a good system of parameters of $M$, and therefore the set of good systems of parameters of $M$ is  non-empty. Let $\x=x_1,\ldots,x_d$ be a good system of parameters of $M$. It easy to see that   $x_1,\ldots,x_{d_i}$ is a good system of parameters of $D_i$ and $x_1^{n_1},\ldots,x_d^{n_d}$ is a good system of parameters of $M$ for any $d$-tuple of positive integers $n_1,\ldots , n_d$. 
\begin{lem}\label{24}
Let $\x=x_1,\ldots,x_d$ be a good system of parameters of $M$. Then $D_i=0\:x_j$ for all  $j= d_i+1, \ldots , d_{i+1}$, $i=0,1,\ldots,t-1$, and therefore
$0\:x_1\subseteq0\:x_2\subseteq\ldots\subseteq0\:x_d.$
\end{lem}
\begin{proof}
Since $D_i\cap(x_{d_i+1},\ldots,x_d)M=0$, we have $D_i\subseteq 0\:x_j$ for all $j\geqslant d_i+1$. Thus it suffices to prove that $0\:x_j\subseteq D_j$ for any $d_i<j\leqslant d_{i+1}$. Assume that $0\:x_j\not\subseteq D_i$. Let $s$ be the largest integer such that $0\:x_j\not\subseteq D_{s-1}$. Then $t\geqslant s>i$ and $0\:x_j=0:_{D_s}x_j$. Since $d_s\geqslant d_{i+1}\geqslant j$, $x_j$ is a parameters element of $D_s$ and therefore $\dim 0\:x_j<d_s$. Hence $0\: x_j\subseteq D_{s-1}$ by the maximality of $D_{s-1}$. This contradicts to the choose of $s$. Therefore $0\:x_j=D_i$. 
\end{proof}
\begin{lem}\label{25}
Let $N$ be a submodule of $M$ such that $\dim N<\dim M$ and $M/N$ a Cohen-Macaulay module. Let $x_1,\ldots,x_i$ be a part of a system of parameters of $M$. Then $$(x_1,\ldots,x_i)M\bigcap N=(x_1,\ldots,x_i)N.$$
\end{lem}
\begin{proof}
We argue by the induction on $i$. The case $i=1$ is trivial. Assume that $i> 1$. Let $a\in (x_1,\ldots,x_i)M\bigcap N$. We write $a=x_1a_1+\ldots+x_ia_i,$ where $a_j\in M,$ $j=1,\ldots,i$. Since $a\in N$, $a_i\in(N+(x_1,\ldots,x_{i-1})M):x_i$. On the other hand, since the sequence $x_1,\ldots,x_i$ is $M/N$-regular, $(N+(x_1,\ldots,x_{i-1})M):_Mx_i=N+(x_1,\ldots,x_{i-1})M$, and  we get $a_i\in N+(x_1,\ldots,x_{i-1})M$. Write $a_i=x_1b_1+\ldots+x_{i-1}b_{i-1}+c,$ where $b_j\in M,$ $j=1,\ldots,i-1$ and $c\in N$. Then  $a-x_ic\in (x_1,\ldots,x_{i-1})M\cap N=(x_1,\ldots,x_{i-1})N$ by the inductive hypothesis. Hence $a\in(x_1,\ldots,x_i)N$. 
\end{proof}

Recall that $M$ is said to be a {\it sequentially Cohen-Macaulay module}, if each quotient $D_i/D_{i-1}$ in the dimension filtration of $M$ is Cohen-Macaulay. Note that the notion of sequentially Cohen-Macaulay modules was introduced first by Stanley in  \cite{St} for the graded case, and  was studied for the local case in \cite{Sch}, \cite{CN}.
The following result is an immediate consequence of Lemma \ref{25} and the definition of a good system of parameters.
\begin{cor}\label{26}
Let $\x=x_1,\ldots,x_d$ be a good system of parameters of a sequentially Cohen-Macaulay $M$. Then $(x_1,\ldots,x_d)M\cap D_i=(x_1,\ldots,x_{d_i})D_i$ for all $i=1,\ldots,t-1$.
\end{cor}

\section{Proof of Theorem 1.1}
To prove Theorem 1.1 we need some auxiliary lemmata. The following result is due to Heinzer-Ratliff-Shah \cite[Theorem 2.4]{HRS}. But we give here the module version of this result proved by  Goto-Shimoda \cite[Lemma 2.1]{GS2}.

\begin{lem}\label{30}
Let $s$ be a positive integer and  $y_1,\ldots,y_s$ an $M$-regular sequence in $\m$. Then
$$(y_1,\ldots,y_s)^nM=\pds$$
for all $n\geqslant 1$.
\end{lem} 

With the same methods that used in \cite{GS1}, we can prove the following results which are module versions of Proposition 3.4 and Lemma 2.1 of \cite{GS1}. 

\begin{lem}\label{28} Let $s$ be a positive integer and $y_1,\ldots,y_s$  a sequence of elements in $\m$ such that $(y_1,\ldots,y_s)^nM=\pds$ for all $n\geqslant 1$.  Then 

\noindent (i) $(y_1,\ldots,y_i)^nM=\pdm$ for all $n\geqslant 1$ and $ i<s$.

\noindent (ii) $y_{i+1}^k M\cap(y_1,\ldots,y_i)^mM\subseteq(y_1,\ldots,y_i,y_{i+1})^{k+m}M$ for all $k, m\geqslant 1$ and $i<s$.
\end{lem}
\begin{lem}\label{27}
Let $s$ be a positive integer and $y_1,\ldots,y_s$ a sequence of elements in $\m$ such that $(y_1,\ldots,y_s)^nM=\pds$ for all $n\geqslant 1$.   Then  
$$y_{i+1}^k M\cap(y_1,\ldots,y_i)^mM\subseteq y_{i+1}^k(y_1,\ldots,y_i,y_{i+1})M+(y_1,\ldots,y_i)^{m+1}M$$ for all $k, m\geqslant 1$ and $1\leqslant i< s$.
\end{lem}
\begin{proof} By Lemma \ref{28}  it is enough  to show that $$(y_1,\ldots,y_i,y_{i+1})^{k+m}M\subseteq y_{i+1}^k(y_1,\ldots,y_i,y_{i+1})M+(y_1,\ldots,y_i)^{m+1}M$$for all $k, m\geqslant 1$.
Indeed, let $a$ be an  element of $ M$ and $(n_1,\ldots,n_i,n_{i+1})\in\Bbb Z^{i+1}$ such that $n_1+\ldots+n_{i+1}=k+m$. If $n_{i+1}\geqslant k$, then $n_1+\ldots+n_i+(n_{i+1}-k)=m\geqslant 1$. Thus 
$$y_1^{n_1}\ldots y_i^{n_i}y_{i+1}^{n_{i+1}}a=y_{i+1}^k(y_1^{n_1}\ldots y_i^{n_i}y_{i+1}^{n_{i+1}-k})a\in y_{i+1}^k(y_1,\ldots,y_i,y_{i+1})M.$$ If $n_{i+1}\leqslant k-1$, then $n_1+\ldots+n_i=k+m-n_{i+1}\geqslant m+1$. Thus $y_1^{n_1}\ldots y_i^{n_i}y_{i+1}^{n_{i+1}}a\in (y_1,\ldots,y_i)^{m+1}M$. Therefore  
$$y_1^{n_1}\ldots y_i^{n_i}y_{i+1}^{n_{i+1}}a\in y_{i+1}^k(y_1,\ldots,y_{i+1})M+(y_1,\ldots,y_i)^{m+1}M,$$ and  the inclusion follows.
\end{proof} 
\begin{lem}\label{13}
Let $\x=x_1,\ldots,x_d$ be a system of parameters of $M$ having the property of parametric decomposition. Then for all $1\leqslant i< j\leqslant d$, there exists an integer $k\geqslant 1$ such that $\q_iM:x_j^n=\q_iM+0\:x_j^k$ for all $n\geqslant k$.
\end{lem}
\begin{proof} First, we claim that $x_j^n M\cap\q_iM\subseteq x_j^n(x_j,\q_i) M$ for all $n\geqslant 1$. Assume the contrary. Then, by
Krull's Intersection Theorem there is an integer $m\geqslant 1$ so that $x_j^n M\cap\q_iM\subseteq x_j^n(x_j,\q_i)M+\q_i^mM$ but $x_j^nM\cap\q_iM\not\subseteq x_j^n(x_j,\q_i)M+\q_i^{m+1}M$. Therefore 
$$x_j^nM\cap\q_iM\subseteq x_j^nM\cap [x_j^n(x_j,\q_i)M+\q_i^mM]=x_j^n(x_j,\q_i)M+x_j^nM\cap\q_i^mM.$$
On the other hand, by Lemma \ref{27} and the hypothesis, we get $x_j^n M\cap\q_i^mM\subseteq x_j^n(x_j,\q_i)M+\q_i^{m+1}M$. It follows that $x_j^n M\cap\q_iM\subseteq x_j^n(x_j,\q)M+\q_i^{m+1}M$, which is impossible. Hence $x_j^n M\cap\q_iM\subseteq x_j^n(x_j,\q_i)M$ and the claim is proved. Thus
$$x_j^n(\q_iM:x_j^n)\subseteq x_j^n M\cap\q_iM\subseteq x_j^n(x_j,\q_i)M.$$
Therefore $\q_iM:x_j^n\subseteq(x_j,\q_i)M+0\:x_j^n$. Take $k\gg0$ so that $\q_iM:x_j^k=\q_iM:x_j^{k+1}$ and $0\:x_j^k=0\:x_j^{k+1}$. Then $\q_iM:x_j^n\subseteq(x_j,\q_i)M+0\:x_j^k$ for all $n\geqslant k$. Let $a\in\q_iM:x_j^n$, we write $a=x_jb+x_1b_1+\ldots+x_ib_i+c$, where $c\in0\:x_j^k$. Since $x_j^na\in\q_iM$ and $n\geqslant k$,  $b\in\q_iM:x_j^{n+1}$. Thus $a\in x_j(\q_iM:x_j^{n+1})+\q_iM+0\:x_j^k$. It follows that $$\q_iM:x_j^n=x_j(\q_iM:x_j^{n+1})+\q_iM+0\:x_j^k =x_j(\q_iM:x_j^n)+\q_iM+0\:x_j^k$$ for all $n\geqslant k$. Hence $\q_iM:x_j^n=\q_iM+0\:x_j^k$ by Nakayama Lemma. 
\end{proof}

 Now we are able to prove  Theorem 1.1.
\begin{proof}[Proof of Theorem 1.1]
$(i)\Rightarrow(ii)$. Let $\x= x_1,\ldots,x_d$ be a good system of parameters of $M$. We  prove by the induction on the length $t$ of the dimension filtration $\mathcal{D}$ of $M$ that $\x$ has the property of parametric decomposition. The case $t=0$ is trivial. Set $\o M=M/D_{t-1}$. Since $\o M$ is a Cohen-Macaulay module, the sequence $x_1,\ldots,x_d$ is $\o M$-regular. Then $\q^n \o M=\bigcap\limits_{\alpha\in\Lambda_{d,n}}\qa\o M$, therefore  $\pd\subseteq \q^nM+D_{t-1}$. Since $x_1^{\alpha_1},\ldots,x_d^{\alpha_d}$ is a good system of parameters of $M$ for all $\alpha \in \Lambda _{d,n}$, it follows by Corollary \ref{25} that  $\qa M\cap D_{t-1}=(x_1^{\alpha_1},x_2^{\alpha_2},\ldots,x_{d_{t-1}}^{\alpha_{d_{t-1}}})D_{t-1}$. Thus  
$$\begin{aligned}
\pd&=[\pd]\cap[\q^nM+D_{t-1}]\\
&=\q^nM+\bigcap\limits_{\alpha\in\Lambda_{d,n}}[\qa M\cap D_{t-1}]\\
&=\q^nM+\bigcap\limits_{\alpha\in\Lambda_{d,n}}(x_1^{\alpha_1},x_2^{\alpha_2},\ldots,x_{d_{t-1}}^{\alpha_{d_{t-1}}})D_{t-1} .
\end{aligned}$$
Note that  $(\beta_1,\ldots,\beta_{d_{t-1}},1,\ldots,1)\in \Lambda_{d,n}$ for any
 $(\beta_1,\ldots,\beta_{d_{t-1}})\in \Lambda_{d_{t-1},n}$ and the length of the dimension filtration of  the sequentially Cohen-Macaulay module $D_{t-1}$ is $t-1$. Therefore, by the inductive hypothesis we have 
$$\begin{aligned}
\bigcap\limits_{(\alpha_1,\ldots,\alpha_d)\in\Lambda_{d,n}}(x_1^{\alpha_1},\ldots,x_{d_{t-1}}^{\alpha_{d_{t-1}}})D_{t-1}&\subseteq \bigcap\limits_{(\beta_1,\ldots,\beta_{d_{t-1}})\in\Lambda_{d_{t-1},n}}(x_1^{\beta_1},\ldots,x_{d_{t-1}}^{\beta_{d_{t-1}}})D_{t-1}\\ &
=(x_1,x_2,\ldots,x_{d_{t-1}})^nD_{t-1} \subseteq \q^nM.\end{aligned}$$
Hence
$$\pd=\q^nM$$ as required.

\noindent
$(ii)\Rightarrow(iii)$ is obvious.

\noindent
$(iii)\Rightarrow(i)$. Let $\x= x_1,\ldots,x_d$ be a good system of parameters  of $M$ having the property of parametric decomposition. We show first that $(\q_iM+D_s):x_{i+1}=\q_iM+D_s$ for all $i< d_{s+1}$ and $s=0,\ldots,t-1$. Indeed, there exists by Lemma \ref{13} a positive integer $k$ such that   $\q_iM:x_{i+1}^k=\q_iM+0\:x_{i+1}^k$ and $\q_iM:x_{d_{s+1}}^{k+1}=\q_iM+0\:x_{d_{s+1}}^{k}$. Observe  by Lemma \ref{24} that $0\:x_{i+1}^k\subseteq 0\:x_{d_{s+1}}^{k}$. Then we have \[\begin{aligned}(\q_iM+0\:x_{d_{s+1}}):x_{i+1}^k&\subseteq\q_iM:x_{d_{s+1}}x_{i+1}^k\\&=(\q_iM+0\:x_{i+1}^k):x_{d_{s+1}})\\&\subseteq\q_iM:x_{d_{s+1}}^{k+1}=\q_iM+0\:x_{d_{s+1}}^k.\end{aligned}\] Note by Lemma \ref{24} that $D_s=0\:x_{d_{s+1}}^k$, so we get 
$$(\q_iM+D_s):x_{i+1}^k\subseteq
\q_iM+D_s\subseteq (\q_iM+D_s):x_{i+1}$$ for all $i< d_{s+1}$, and the conclusion follows.
This implies that  $\depth M/D_s\geqslant d_{s+1}$ for $s=0,\ldots,t-1$. Now, from    the short exact sequences
$$0\to D_s/D_{s-1}\to M/D_{s-1}\to M/D_s\to0,$$ 
it follows that $D_s/D_{s-1}$ is Cohen-Macaulay for all $s=1,\ldots,t$, and the proof of Theorem 1.1 is complete. 
\end{proof}

\begin{proof}[Proof of Corollary 1.2]
$(i)\Rightarrow(ii)$. It is easy to see from the hypothesis that $M$ is a sequentially Cohen-Macaulay module with the dimension filtration $\D: H^0_\m(M)\subset M$. Moreover, by Lemma \ref{25} we have $$(x_1,\ldots,x_d)M\cap H^0_\m(M)=(x_1,\ldots,x_d)H^0_\m(M)\subseteq\m H^0_\m(M)=0$$ for any system of parameters $x_1,\ldots,x_d$ of $M$. This means that every system of parameters of $M$ is good, therefore it has the property of parametric decomposition by Theorem 1.1.

\noindent
$(ii)\Rightarrow(i)$. First, it follows by Theorem 1.1 that $M$ is  sequentially Cohen-Macaulay. Remember that  the definition of the dimension filtration of $M$ that $D_0= H^0_\m(M)$ and $\dim D_i >0$ for all $i>0$. Therefore the implication is proved, if  we can show that $\m D_{t-1}=0$. 
Suppose the contrary.  Then there is an element $x_1\in\m$ so that $x_1D_{t-1}\not=0$ and $\dim M/x_1M=d-1$. Since $d\geqslant 2$, we can choose $x_2\in\m$ such that $x_2D_{t-1}=0$ and $\dim M/(x_1,x_2)M=d-2$. We observe that the sequence $x_1,x_2$ and $x_1,x_1+x_2$ are part of systems of parameters of $M$. Therefore, by the hypothesis and Lemma \ref{28}, (i) we get 
 \[\begin{aligned}(x_1^2,x_1+x_2)M\cap(x_1,(x_1+x_2)^2)M&=(x_1,x_1+x_2)^2M\\&=(x_1,x_2)^2M=(x_1^2,x_2)M\cap(x_1,x_2^2)M.\end{aligned}\]
Since $M/D_{t-1}$ is Cohen-Macaulay,  it follows from  Lemma \ref{25} that
\[\begin{aligned} x_1D_{t-1}&=(x_1^2,x_1+x_2)D_{t-1}\cap(x_1,(x_1+x_2)^2)D_{t-1}\\&=(x_1^2,x_2)D_{t-1}\cap(x_1,x_2^2)D_{t-1}=x_1^2D_{t-1}.\end{aligned}\]
 Thus $x_1D_{t-1}=0$ by Nakayama's lemma,  which is impossible. Hence $\m D_{t-1}=0$.  
\end{proof}
\section{Hilbert-Samuel polynomials}
A parameter ideal $\q$ is called a {\it good parameter ideal} if it is generated by a good system of parameters. Then, in this section we shall show  that for a sequentially Cohen-Macaulay module $M$ the Hilbert-Samuel function $H_{\q, M}(n)=\ell(M/\q^{n+1}M)$ has a special expression with non-negative coefficients, which can be computed by the dimension filtration, and this function coincides with the Hilbert-Samuel polynomial  $P_{\q, M}(n)$ for any good parameter ideal $\q$ of $M$ and all $n\geqslant 1$. Moreover, the sequentially Cohen-Macaulayness of $M$ can be characterized by this expression of the Hilbert-Samuel function. First, we begin with the following lemma which is an easy consequence of Theorem 1.1.
\begin{lem}\label{41}
Let $\q$ be a good parameter ideal of a sequentially Cohen-Macaulay module $M$. Then 
$$\q^nM\cap D_i=\q^nD_i,$$
for all $n\geqslant 1$ and $i=0,\ldots,t$.
\end{lem}
\begin{proof}
Since $\q$ is a good parameter ideal of $M$, there is a good system of parameters $x_1,\ldots,x_d$ of $M$ such that $\q M=(x_1,\ldots,x_d)M$. Then by Theorem 1.1 and Corollary \ref{26}, we get \[\begin{aligned}\q^nM\cap D_i&=\big[\pd\big]\cap D_i\\&=\bigcap\limits_{\alpha\in\Lambda_{d,n}}(\qa M\cap D_i)=\bigcap\limits_{\alpha\in\Lambda_{d,n}}(x_1^{\alpha_1},x_2^{\alpha_2},\ldots,x_{d_i}^{\alpha_{d_i}})D_i.\end{aligned}\]  
Observe that $(\beta_1,\ldots,\beta_{d_i},1,\ldots,1)\in \Lambda_{d,n}$ for all $(\beta_1,\ldots,\beta_{d_i})\in \Lambda_{d_i,n}$. Therefore,  we obtain by Theorem 1.1 that 
 \[\begin{aligned}
\bigcap\limits_{\alpha\in\Lambda_{d,n}}(x_1^{\alpha_1},x_2^{\alpha_2},\ldots,x_{d_i}^{\alpha_{d_i}})D_i& \subseteq\bigcap\limits_{(\beta_1,\ldots,\beta_{d_i})\in\Lambda_{d_i,n}}(x_1^{\beta_1},x_2^{\beta_2},\ldots,x_{d_i}^{\beta_{d_i}})D_i\\&=(x_1,\ldots,x_{d_i})^nD_i.
\end{aligned}\]
 So $\q^nM\cap D_i\subseteq(x_1,\ldots,x_{d_i})^nD_i\subseteq\q^nD_i$ and the conclusion follows.
\end{proof}
The following result seems to be well-known. But, as we can not find a reference to it,  we give a brief proof for the sake of completeness.
\begin{lem}\label{42}
Let $\q$ be a parameter ideal of module $M$. Then
$$\ell(M/\q^{n+1}M)\leqslant\binom{n+d}{d}\ell(M/\q M).$$
Moreover, this inequality becomes an equality if and only if $M$ is a Cohen-Macaulay module. 
\end{lem}
\begin{proof}
Suppose that $\q=(x_1,\ldots,x_d)$ is a parameter ideal of $M$. We set $N=(M/\q M)[X_1,\ldots,X_d]$ and $\gr_\q(M)=\bigoplus\limits_{i=0}^{\infty}\q^iM/\q^{i+1}M$. Then one has the natural surjection $\varphi: N\to \gr_\q(M)$ defined  by  $\varphi(X_i)= \bar x_i\in\q/\q^2$. Therefore
$$\ell(M/\q^{n+1}M)\leqslant \ell (N/(X_1,\ldots , X_d)^{n+1}N)=\binom{n+d}{d}\ell(M/\q M).$$
Moreover, the last inequality becomes an equality if and only if $\varphi$ is an isomorphism, and this condition is clear equivalent to the Cohen-Macaulayness of $M$.
\end{proof}
 
\begin{theorem}\label{43} Let $\D: D_0\subset D_1\subset\ldots\subset D_t=M$ be the dimension filtration of $M$ and set $\D_i=D_i/D_{i-1}$ for all $i=1,\ldots,t$, $\D_0 =D_0$. Then
the following statements are equivalent:

\noindent (i) $M$ is a sequentially Cohen-Macaulay module.

\noindent (ii) For any good parameter ideal $\q$ of $M$, it holds $$\ell(M/\q^{n+1}M)=\sum\limits_{i=0}^t\binom{n+d_i}{ d_i}\ell(\D_i/\q\D_i)$$
for all $n\geqslant0$. 

\noindent (ii) There exists  a good parameter ideal $\q$ of $M$ such that
$$\ell(M/\q^{n+1}M)=\sum\limits_{i=0}^t\binom{n+d_i}{ d_i}\ell(\D_i/\q\D_i)$$
for all $n\geqslant0$.   
\end{theorem}

\begin{proof}
 $(i)\Rightarrow(ii)$. We argue by the induction on the length $t$ of the dimension filtration $\mathcal{D}$ of $M$. The case $t=0$ is obvious.  Assume that $t>0$. By virtue of Lemma \ref{41}, we have a short exact sequence
$$0\to D_{t-1}/\q^{n+1}D_{t-1}\to M/\q^{n+1}M\to M/\q^{n+1}M+D_{t-1}\to0.$$
Therefore, we have $\ell(M/\q^{n+1}M)=\ell(D_{t-1}/\q^{n+1}D_{t-1})+\ell(\D_t/\q^{n+1}\D_t)$. Since $D_{t-1}$ is a sequentially Cohen-Macaulay module and its dimension filtration is of the length $t-1$, it follows from the inductive hypothesis that
$$\ell(D_{t-1}/\q^{n+1}D_{t-1})=\sum\limits_{i=0}^{t-1}\binom{n+d_i}{d_i}\ell(\D_i/\q\D_i).$$
Note that $\D_t$ is Cohen-Macaulay of dimension $d=d_t$, we have $$\ell(\D_t/\q^{n+1}\D_t)=\binom{n+d}{d}\ell(\D_t/\q\D_t).$$ Hence$$\ell(M/\q^{n+1}M)=\sum\limits_{i=0}^t\binom{n+d_i}{ d_i}\ell(\D_i/\q\D_i),$$
for all $n\geqslant0$ as required.

\noindent $(ii)\Rightarrow(iii)$ is trivial. 

\noindent $(iii)\Rightarrow(i)$. Since the following sequence is exact
$$D_{t-1}/\q^{n+1}D_{t-1}\to M/\q^{n+1}M\to M/\q^{n+1}M+D_{t-1}\to0,$$
we get $\ell(M/\q^{n+1}M)\leqslant\ell(D_{t-1}/\q^{n+1}D_{t-1})+\ell(\D_{t}/\q^{n+1}\D_t)$. Therefore, by induction on the length of the dimension filtration we can show that $$\ell(M/\q^{n+1}M)\leqslant\sum\limits_{i=0}^t\ell(\D_i/\q^{n+1}\D_i) .$$
On the other hand, since $$\ell(\D_i/\q^{n+1}\D_i)\leqslant\binom{n+d_i}{ d_i}\ell(\D_i/\q\D_i)$$ for all $i=0,\ldots,t$ by Lemma \ref{42},
 $$\ell(M/\q^{n+1}M)\leqslant\sum\limits_{i=0}^t\ell(\D_i/\q^{n+1}\D_i) \leqslant\sum\limits_{i=0}^t\binom{n+d_i}{ d_i}\ell(\D_i/\q\D_i).$$
It follows from the hypothesis of (iii) that $\ell(\D_i/\q^{n+1}\D_i)=\binom{n+d_i}{ d_i}\ell(\D_i/\q\D_i)$ for all $i=0,\ldots,t$. Thus  $\D_i$ is  Cohen-Macaulay  for all $i=0,\ldots,t$ by Lemma \ref{42} again, and this completes the proof.
\end{proof}

\end {document}